\documentclass[a4paper]{article}

\usepackage{amssymb}
\usepackage{amsmath}
\usepackage{amsfonts,amsthm}
\newtheorem{theorem}{Theorem}

\newtheorem{conjecture}{Conjecture}
\newtheorem{corollary}{Corollary}

\newtheorem{definition}{Definition}

\newtheorem{lemma}{Lemma}

\newtheorem{proposition}{Proposition}

\begin{document}

\title{Strange images of profinite groups}
\author{Nikolay Nikolov}
\date{}
\maketitle

\begin{abstract}
We investigate whether a finitely generated profinite group $G$ could have a finitely generated infinite image. A result of Dan Segal shows that this
is impossible if $G$ is prosoluble. We prove that such an image does not
exist if $G$ is semisimple or \emph{nonuniversal}. We also investigate the
existense of dense normal subgroups in $G$.
\end{abstract}

\section{Introduction} What can the images of a profinite group $G$ be? There
are obviously all the continuous images of $G$ but there may well be others.
In the first place $G$ may not be finitely generated and then some of its
images may be very different from the continuous ones. For example there are examples of profinite
groups which are inverse limits of finite perfect groups while $G \not =
G^2$,
see \cite{danbook}. However this cannot happen when $G$ is finitely generated: \cite{NS} proves
that some verbal subgroups of $G$ are automatically closed: for example any term of the lower central series of $G$. However not all
verbal subgroups of $G$ are closed - the second derived subgroup $G''$
may not be closed, even when $G$
is a finitely generated pro-$p$ group. See \cite{andrei} for a description
of the words $w$ such that $w(G)$ is closed in all finitely generated pro-$p$
groups $G$. 

Let us call a quotient $G/N$ of a profinite group $G$ \emph{strange} if $N$ is not closed in $G$. The strange images of $G$ are defined to be groups
$I$ such that there exists an abstract epimorphism $f: G \rightarrow I$ but
there is no continuous such choice for $f$. Clearly every strange image of
$G$ is also a strange quotient but the converse need not hold (for example
in $G=C_p^N$).

We are motivated by the following
\begin{theorem}[Nikolov \& Segal, \cite{NS}] \label{ns} A finitely generated profinite group has no strange finite quotients. 
\end{theorem}

\begin{corollary} \label{cor1} A residually finite image of a finitely generated profinite
group $G$ is either finite or uncountable.
\end{corollary}

Indeed, if $f: G \rightarrow I$ is a homomorphism with residually finite
image $I$ then $\ker f$ is an intersection of subgroups of finite index in
$G$ and therefore by Theorem \ref{ns} is a closed subgroup of $G$. Hence $I$ is a profinite group
and must be uncountable if it is infinite (say by the Baire category theorem). 
\medskip

It is natural to try to extend this result to other types of images. As
a first step we could ask whether a finitely generated profinite group can have a countable infinite
image. This may happen, in fact already $\mathbb Z_p$ has
such an image: compose the inclusion $\mathbb Z_p <\mathbb Q_p$ with a surjective
linear map $\mathbb Q_p \rightarrow \mathbb Q$ as a vector space over the rationals.

The following conjecture arose in Blaubeuren in 2007 in discussions with Y.
Barnea, E. Breullard, P.E. Caprace, T. Gelander and J. Wilson.

\begin{conjecture}[Blaubeuren] \label{b} A profinite group $G$ cannot have an infinite finitely generated image.
\end{conjecture}

Suppose that the homomorphisms $f: G \rightarrow I$ has finitely generated infinite image. Without loss of generality we may assume that $G$ is finitely
generated as a profinite group. Let $J$ be the intersection of all subgroups
of finite index in $I$. Since $I/J$ is a residually finite image of $G$ by Corollary \ref{cor1} we deduce that $I/J$ is finite. Then $H:=f^{-1}(J)$ is a subgroup of finite index in $G$ hence by Theorem \ref{ns} it is a finitely generated profinite group and $J$ is a finitely generated infinite image of $H$. Now an easy application
of Zorn's lemma shows that $J$, hence $H$ has an infinite simple image.

This shows that it is enough to prove the Conjecture \ref{b} for a finitely generated \emph{simple} image.  

First we consider prosoluble groups. There a theorem of Dan Segal \cite{dan}
gives the answer.

\begin{theorem}[D. Segal] \label{sol}Suppose $G$ is a prosoluble group generated by $d$ elements.
Let $a_1, a_2, \ldots ,a_r$ generate $G$ modulo $G'$, i.e. $G=\overline{\langle a_1, a_2 \ldots ,a_r \rangle} G'$ (note that $G'$ is closed in $G$ by Corollary
1 of \cite{dan}). Then
\[ G'=  \left(\prod_{i=1}^r [G,a_{i}] \right)^{*(72d+46)},
\]
where $[G,y]=\{g^{-1}y^{-1}gy \ | \ g \in G\}$ and for a subset $A \subset
G$ and integer $t$ $A^{*t}:=A \cdot A \cdots A$ ($t$ times). 
\end{theorem}

So if $G/N$ is a simple infinite quotient of $G$ we have that $G'N=G$ and hence we can
choose elements $n_1, \ldots n_d$ from $N$ which generate $G$ modulo $G'$.
Now Theorem \ref{sol} gives that every element of $G'$ is a product of $2d(72d+46)$ conjugates of the elements $n_i$, hence $G' \leq N$, i.e. $N=G$, contradiction.

Therefore Conjecture \ref{b} holds for prosoluble groups $G$. \bigskip

Another class of groups where we can test the validity of Conjecture \ref{b}
is the class of semisimple profinite groups, i.e. the Cartesian products of nonabelian
finite simple groups. For these groups we have the following result.

\begin{proposition} \label{main}
Every simple image of a Cartesian product $G$ of nonabelian finite simple groups
is either finite or uncountable. In particular Conjecture \ref{b} holds for
$G$.
\end{proposition}

In fact in Proposition \ref{main2} below we shall describe all the simple quotients of $G$. 
\bigskip

A natural generalization of prosoluble groups is the class of \emph{nonuniversal} profinite groups.

For a finite group $\Gamma$ let $\alpha(\Gamma)$ be the largest integer $k$
such that $\Gamma$ involves the alternating group $A_k$, i.e. $A_k$ is an image of a subgroup of $\Gamma$.
We say that a profinite group $G$ is \emph{non-universal} if there is an integer
$k$ such that $\alpha (\bar G) \leq k$ for all finite continuous images $\bar
G$ of $G$. In that case we denote by $\alpha(G)$ the maximum of the numbers
$\alpha(\bar G)$. It is easy to see that a profinite group $G$ is non-universal if and only if its continuous finite images don't involve some finite group.

\begin{theorem} \label{nou} Conjecture \ref{b} holds for nonuniversal profinite groups $G$.
\end{theorem}

Proposition \ref{main} together with Theorems \ref{sol} and \ref{nou}
gives the following 
\begin{corollary} Let $G$ be a profinite group which has a subnormal chain $G=G_1 \vartriangleright
G_2 \vartriangleright \cdots \vartriangleright G_k=1$ of closed subgroups
$G_i$ such that each quotient $G_i/G_{i+1}$ is either semisimple
or non-universal.
Then $G$ cannot have a finitely generated infinite image.
\end{corollary}
\subsection{Dense normal subgroups} \label{den}
Notice that when $G/N$ is a strange simple quotient of a profinite group $G$
then $N$ must be a dense subgroup of $G$.

We now turn to  the more general task of investigating dense normal subgroups of profinite groups. This is not such a strange topic as it might appear at first and is of some interest to the study of topologically simple compactly
generated locally compact groups, see section \ref{lc} below.

Of course hoping for a  complete classification of these is too much: Any infinite abelian profinite group
contains many dense normal subgroups. As another example consider a Cartesian product $\prod_{i=1}^\infty
S_i$ of profinite or finite groups $S_i$ which contains the direct product $\bigoplus_{i=1}^\infty
S_i$. So if a profinite group $G$ has an infinite image which is abelian or a Cartesian product of finite groups then $G$ will possess a proper
dense normal subgroup. It is surprising that the converse of this holds,
at least in the case when $G$ is nonuniversal.

\begin{theorem} \label{exdense}
Let $G$ be a finitely generated nonuniversal profinite group.
Then $G$ has a dense proper normal subgroup if and only if one of the following
holds: 

(a) $G/G'$ is infinite or,

(b) $G$ maps onto infinitely many finite simple groups.
\end{theorem}  
 
One might wonder if the condition that $G$ is nonuniversal above can be removed.
The following example shows that it is necessary at least if one wants to
keep the same conclusion of the theorem. 

Consider the semidirect product
$G = \prod_{n=5}^\infty {A_n} \rtimes C_2$ where the generator
$\tau$ of $C_2$ acts on each of the factors $A_n$ as conjugation by a transposition.
Then $G/G' \simeq C_2$ and the only finite simple image of $G$ is $C_2$.
However it is easy to see that the normal closure of $\tau$ in $G$ is a
proper dense subgroup. 

This shows that when $G$ is universal some extra possibilities should be added to the list in Theorem \ref{exdense}. As a first step one might wish to study this in the case when $G$ is virtually semisimple.   

\begin{definition}
Let $G_*$ be the intersection of all open normal subgroups $K$ of $G$ such
that $G/K$ is abelian or simple. 
\end{definition}

Thus $G/G_*$ is the largest quotient of
$G$ which has the form $A \times \prod_{i} S_i$ where $A$ is an
abelian profinite group and $\{S_i\}_i$ are nonabelian simple groups. Theorem \ref{exdense} seems to suggest that any normal dense normal subgroup of $G$ contains $G_*$.
However this is false even for nilpotent groups: 

Let 
$G=U_3( \widehat{\mathbb Z})$ be the Heizenberg group of upper 3-by-3 unitriangular matrices with
entries in $\widehat{\mathbb Z}$. Then 
$G \simeq \prod_{p \textrm{ prime}}U_3(\mathbb{Z}_p)$ and just as before this contains the dense normal subgroup $N=\bigoplus_{p \textrm{ prime}}U_3(\mathbb Z_p)$. Since $N$ does not contain $G'$ the quotient
$G/N$ is not abelian. Similar examples show that $N$ does not need to contain
$G_*$ in some perfect profinite groups $G$ as well.

However in the special case when $N$ contains a finitely generated
normal subgroup then indeed $G_* \leq N$, in fact $G' \leq N$.

\begin{theorem} \label{g'}Let $G$ be a finitely generated nonuniversal profinite group and let $N$
be a normal subgroup of $G$. Suppose that $N$ contains a finitely generated
subgroup $D$ whose image is dense in $G/G_*$. Then $G/N$ is abelian.  
\end{theorem}

\begin{corollary} \label{finitelymany}Let $G$ be a finitely generated nonuniversal
profinite group which
maps onto only finitely many finite simple groups (abelian or nonabelian).
Then any dense normal subgroup of $G$ contains $G'$.
\end{corollary}

\textbf{Proof of Corollary \ref{finitelymany}:} The assumption implies that
the Frattini subgroup $F/G_*$ of  $G/G_*$ is open in $G/G_*$, i.e. $F$ is
open in $G$. Hence if $N$ is a dense normal subgroup of $G$ then $FN=G$.
Therefore we can find finitely many elements of $N$ which generate $G/F$ and hence their images generate a dense subgroup of $G/G_*$. Therefore Theorem
\ref{g'} applies. $\square$

\subsection{Topologically simple locally compact groups}\label{lc}
As promised we give an application of Corollary \ref{finitelymany}.
A locally compact group $G$ is \emph{locally finitely generated} if its compact
open subgroup $U$ is topologically finitely generated, i.e. contains a dense
finitely generated subgroup.

\begin{proposition} \label{top} Let $G$ be a compactly generated totally disconnected locally compact topological
group which is locally finitely generated. Assume that $G$ is topologically
simple. Then the only proper abstract quotients of $G$ must be abelian.
\end{proposition} 

Proposition \ref{top} can be viewed as an analogue of the following result
in the connected case (cf. Theorem XVI 2.1 in \cite{Ho}): If $N$ is a dense analytic subgroup
of Lie group $G$ there is an abelian subgroup $A$ of $G$ such that
$G=NA$. \bigskip

\textbf{Proof of Proposition \ref{top}.}
Let $U$ be a compact open subgroup of $G$. Then $U$ is a profinite group
which is finitely generated. Let $K$ be a compact generating set for $G$.
Without loss of generality we may assume that $K=K^{-1}=UKU$.
As explained on page 150 of \cite{Mo} $G$ acts on the regular Schreirer graph
$\Gamma$ with vertices $G/U$ and edges $(gU,gkU)$ for $g \in G, k \in K$. Moreover, since $G$ is topologically simple the action is faithful.

Now the
valency of every vertex of $\Gamma$ is $r=|KU/ U|$ which is finite. The profinite
group $U$ is the stabilizer of a vertex of $\Gamma$ and is therefore a quotient of
the group of rooted automorphisms of the infinite $r$-regular tree (the universal
cover of $\Gamma$). It follows that every upper composition factor of $U$
is a quotient of a subgroup of $\mathrm{Sym}(r)$, in particular $U$ is nonuniversal and
maps onto finitely many finite simple groups. Now suppose that $N$ is a nontrivial
normal subgroup of $G$. Since $G$ is topologically simple $N$ must be dense
in $G$ and hence $N_1=N \cap U$ is a dense normal subgroup in $U$. By Corollary
\ref{finitelymany} it follows that $G/N_1$ must be abelian. Finally note that $NU=G$ since $N$
is dense in $G$ and therefore $G/N = NU/N \simeq U/N_1$. $\square$

I wish to thank Pierre-Emmanuel Caprace \cite{PE} for pointing out this application
to me. In light of this result it will be interesting to know
if it is possible that the groups $G$ above can ever have nontrivial abelianization.

The structure of the rest of the article is as follows. Proposition \ref{main} is proved in Section 1 and the proof of Theorem \ref{nou} is in Section \ref{sec2}. Theorems \ref{exdense} and \ref{g'} are proved
in Section \ref{sec3}. 
First we need a few definitions.

\subsubsection*{Notation:}
For an element $g$ in a group $G$ the conjugacy class of $g$ is denoted $g^G$.
When $A,B \subset G$ by $AB$ we denote the set $\{ab \ | \ a
\in A, b \in B \}$ and for $t \in \mathbb N$ let
\[A^{*t}= A \cdot A  \cdots \cdot A \quad (t \textrm{ times }).\]

If $f \in \mathrm{Aut}(G)$ then $[G,f]:=\{g^{-1} g^f
\ | \ g \in G\}$. In particular if $N \triangleleft G$ and $a \in G$ then
$[N,a]:=\{n^{-1} n^a \ | \ n \in N\}$.

Further let $[N,_n G]$ be defined by 
\[ [N,_1 G]=[N,G]= \langle \{[N,g] \ | \ g \in G \}\rangle , \quad [N,_{n+1} G]=[[N,_n
G],G] , \ (n \geq 1)\] and finally let 
\[ [N,_\omega G]= \bigcap_{i=1}^\infty [N,_n G]. \]

We call a (finite or profinite) group \emph{semisimple} if it is a Cartesian product of finite simple groups.
A finite group $G$ is quasisimple if $G=G'$ and $G/\mathrm{Z}(G)$ is simple.

A group $G$ is a \emph{central product} of its subgroups $\{S_i\}_{i \in I}$ if $G= \langle \{S_i\}_i \rangle $ and $[S_i,S_j]=1$ for all distinct $i,j \in
I$. Then $G / \mathrm{Z}(G)$ is isomorphic to the cartesian product of $S_j
/ \mathrm{Z}(S_j)$. A central product of quasisimple groups is called quasi-semisimple.
In that case $\mathrm{Aut}(G)$ acts faithfully on $G/\mathrm{Z}(G)$.

\section{Proof of Proposition \ref{main}}
Let $\{S_i\}_{i=1}^\infty$ be a collection of nonabelian finite simple groups
and let $G= \prod_{i} S_i$. Take an ultrafilter $\mathcal U$ on the integers and recall
the definition of the \emph{ultralimit}, $\lim_\mathcal{U} a_i$ of a 
bounded sequence $(a_i)_{i=1}^\infty$
of real numbers: \ $\lim_\mathcal{U} a_i$ is the unique number $\alpha$ such
that for any $\epsilon >0$ the set of indices $\{ i \in \mathbb N \ | \ |a_i-\alpha|<\epsilon
\}$ is a member of $\mathcal U$. \bigskip

Define a function $h=h_\mathcal U : G \rightarrow [0,1]$ by 
\[h (\mathbf g)= \lim_{\mathcal U} \frac{\log |g_i^{S_i}|}{\log |S_i|}  \ \textrm{ for any } \ \mathbf g= (g_i)_i \in
G=\prod_{i=1}^\infty S_i.\]

It is easy to see that $h(ab) \leq h(a) + h(b)$ for any $a, b \in G$ and
therefore the set 
\[K=K_{\mathcal U}:=h^{-1}(0)\] is a normal subgroup of $G$.
  
\begin{proposition} \label{main2} The simple images of $G = \prod_{i=1}^\infty
S_i$ are exactly the quotients $G/K_\mathcal{U}$ for some ultrafilter $\mathcal
U$ on the integers.
\end{proposition}

 \medskip
 
\textbf{Proof of Proposition \ref{main2}}
As a first step we show that each quotient $G/K_\mathcal{U}$ is indeed simple.
Suppose that $\mathbf g= (g_i)_i \in G$ is outside $K=K_{\mathcal U}$. This means that $h(\mathbf
g)>0$, i.e. the set 
\[A=\{i \in \mathbb N \ | \ \log |g_i^{S_i}|/\log |S_i|>e\}\] is
in $\mathcal{U}$ for some $e>0$.

Now we use the following result 

\begin{theorem}[Liebeck \& Shalev, \cite{LS}]\label{lie} There is a constant $c>0$ such that if
$S$ is a nonabelian finite simple group and $\mathcal C \subset S$ is a  conjugacy class in $S$ then 
\[ S= \mathcal C^{*n},\]
where $n= [c \frac{\log |S|}{\log |\mathcal C|}]$.
\end{theorem}

With this result by setting $n=[c/e]$ we obtain that $S_i= g_i^{S_i} \cdot g_i^{S_i} \cdots g_i^{S_i}$ ($n$ times) for each $i \in A$,
and therefore  
\[G= K \cdot \left( \mathbf{g}^{G} \right)^{*n} \]
because $A \in \mathcal U$.

Since $\mathbf g$ was an arbitrary element from $G$ outside $K$ we see
that $G/K$ is simple. \bigskip

Conversely, suppose that $G/H$ is simple for a proper normal subgroup $H$
of $G$.
For every $\mathbf{t}=(t_i)_i \in H$ and every $\epsilon>0$ let 

\[A(\mathbf t, \epsilon)=\left \{ i \in \mathbb N \ | \ \frac{\log |t_i^{S_i}|}{\log
|S_i|}< \epsilon \right \} \subseteq \mathbb N \] and let $U$ be the collection of all subsets $A(\mathbf t, \epsilon)$ for all $\mathbf t \in H$ and all $\epsilon
>0$.

I \textbf{claim} that every finite subset of $U$ has nontrivial intersection.
Indeed, suppose that 
\[ A(\mathbf {t}_1, \epsilon_1) \cap A(\mathbf {t}_2, \epsilon_2) \cap \cdots
\cap A(\mathbf {t}_k, \epsilon_k)= \emptyset \] for some $A(\mathbf t_i , \epsilon_i)
\in U$. Let $\epsilon= \min_i \{\epsilon_i\}$ and write $\mathbf t_i=(t_{i,j})_j$
with $t_{i,j} \in S_j$.
  
It follows that for each integer $j \in \mathbb N$ there is some $i \leq
k$ such that $j \not \in A(\mathbf t_i, \epsilon)$, i.e., $\log |t_{i,_j}^{S_j}|/\log|
S_j|\geq\epsilon $. Put $N= c/\epsilon$, then Theorem \ref{lie} gives 
\[ S_j=\left( t_{i,j}^{S_j} \right)^{*N}.
\]

By considering independently each coordinate
$j\in \mathbb N$ we see that

\[G= \prod_{i=1}^k \left( \mathbf t_i^G\right)^{*N}.\]

However the elements $\mathbf t_i$ are from $H \vartriangleleft G$ and hence $H=G$. This contradiction
proves the claim.

Since the finite collections of elements of $U$ have nontrivial intersection
a standard application of Zorn's lemma provides the existence of an ultrafilter
$\mathcal U$ on $\mathbb N$ containing $U$. Now from the definition of $U$
it follows that $h_\mathcal{U}(\mathbf t)=0$ for all $\mathbf t \in H$, i.e., $H \leq K_\mathcal
U$. But $G/H$ is simple and therefore $H=K_{\mathcal U}$. Proposition \ref{main2}
is proved.

\bigskip

\textbf{Proof of Proposition \ref{main}}
Armed with the description of the simple images of $G$ it is easy to complete
the proof.
\medskip

In the rest of this section the rank of a finite simple group
$S$ is its (untwisted) Lie rank if $S$ has Lie type, it is $n$ if $S=A_n$ and $0$ otherwise. 

There are two cases to consider.
\medskip

\textbf{Case 1.} For some $m \in \mathbb N$ the collection of indices
$D(m):=\{ i\ | \ \mathrm{rank}(S_i) \leq m \}$ is in $\mathcal U$.

Since ultrafilter $\mathcal U$ is supported on the set $D(m)$, we can assume that $D(m)=\mathbb N$. Now a nontrivial conjugacy class $\mathcal C$ of
a finite simple group $S_i$ has size at least $|S_i|^\epsilon$ for some $\epsilon$ which depends
only on the rank of $S_i$. Therefore the quotient $G/K_\mathcal U$ in this case
coincides with the ultraproduct $\prod_i S_i/\mathcal U$ and by \cite{ultra} an ultraproduct
of finite sets is either finite or uncountable. So it is not possible for
the image to be countably infinite.
\medskip

\textbf{Case 2.} $D(m) \not \in \mathcal U$ for all $m \in \mathbb N$. 

Note that the map $h: G \rightarrow [0,1]$ induces a well defined map on
$G/K_{\mathcal U}$. I claim that the image of $h$ is the whole interval $[0,1]$
and so $G/K_\mathcal U$ must be uncountable.
In turn $h(G)=[0,1]$ follows from
 
\begin{proposition}\label{image} Let $S_i$ be a collection of finite simple
groups such that $\mathrm{rank}(S_i) \rightarrow \infty$ as $i \rightarrow \infty$. For any real number $\beta \in [0,1]$ there is a sequence of elements $g_i \in S_i$ such that $\frac{\log |g_i^{S_i}|}{\log |S_i|} \rightarrow
\beta$ as $i \rightarrow \infty$.

\end{proposition}

Assuming this for the moment let the elements $g_i$ be as in Proposition \ref{image} and let $\mathbf g= (g_i)_i \in G$. Then $h_\mathcal U (\mathbf g)= \beta$ for any nonprincipal ultrafilter
$\mathcal U$. Case 2 is completed and so is the proof of Proposition \ref{main2}. $\square$ \bigskip

\textbf{Proof of Proposition \ref{image}.} 

Put $\alpha= 1- \beta$.

If $S_i$ is an alternating group
$A_n$ take $g_i$ to be an even cycle of length about $\alpha n$ in $A_n$.
\medskip

If $S=S_i$ is a simple classical group, consider the corresponding quasisimple
classical group $\tilde S$ acting on its natural module $V$ over a finite field
of size $q$ equipped with
a bilinear form $f$ (symmetric, sesquilinear, alternating or just equal to 0 in case $S_i$ has type $PSL_n$). We have $S = \tilde S /Z$ where $Z$ is the
centre of $\tilde S$ and if
$g= \tilde g Z \in S$ with $\tilde g \in \tilde S$ we have that

\[ |g^S| \leq |\tilde g^{\tilde S}| \leq |Z| |g^S|\] 

Since the centre $Z$ of $S$ has asymptotically negligible size compared to the size of $S$
it is enough to find an element $\tilde g \in \tilde S$ with 
$\log |\tilde g ^{\tilde S}| \sim \beta \log |\tilde S|$.  
\bigskip

We can decompose $V$ as $V_0 \oplus V_1
\oplus V_2$ such that:

\begin{itemize} 
\item $\dim V_0$ is about $\sqrt{\alpha} \dim V$, and $\dim V_1 =
\dim V_2$, 
\item $V_1 \oplus
V_2$ is orthogonal to $V_0$, and 
\item The form $f$ is nondegenerate on both $V_0$ and $V_1 \oplus V_2$ and is isotropic on $V_1$ and on $V_2$

\end{itemize}

Let $\tilde g \in \tilde S$ be equal to the identity on $V_0$ and act on each of $V_1$
and $V_2$ as a cyclic transformation without fixed vectors. In other words
there is a vector $v_i \in V_i$, $(i=1,2)$ such that $v_i, \tilde gv_i, \tilde
g^2v_i, \ldots $ is a basis for $V_i$.

Now $C_{\tilde S}(\tilde g)$ contains the classical group $L$ on $V_0$ preserving $f$ and
by the choice of $\dim V_0$ we have $\log |L|/\log |\tilde S| \sim  (\dim V_0/\dim
V)^2$ which tends to $\alpha$ as $\dim V \rightarrow \infty$.

On the other hand if $s \in \tilde{S}$ commutes with $\tilde g$ then $s$ must stabilize
$V_0$, the fixed space of $\tilde{g}$. Since $V_1$ and $V_2$ are cyclic modules for
$\tilde{g}$ we have that the action of $s$ on $V_1$ and $V_2$ is determined from
$s\cdot v_1$ and $s\cdot v_2$. Hence $s$ is completely known from its restriction
to $V_0$ and from the two vectors $sv_1,sv_2 \in V$. Denote by $\mathrm{Gf}(V_0)$
the subgroup of $\mathrm{GL}(V_0)$ which preserves $f$. We have that 
$|\mathrm{Gf}(V_0)| \leq q |L|$.

Therefore \[ |L| \leq C_{\tilde S}(\tilde g) \leq |\mathrm{Gf}(V_0)||V|^2 \leq q^{1+ 2 \dim V} |L|
\]
which gives $\log |G_{\tilde S}(\tilde g)|/ \log |\tilde S| \sim \log |L|
/\log |\tilde S| \rightarrow \alpha$ as $\dim V$ tend to infinity.

This easily gives $\log |g^S|/\log|S| \rightarrow 1-\alpha=\beta$ as the
rank of $S$ tends to infinity and proves Proposition \ref{image}.   

\section{Proof of Theorem \ref{nou}} \label{sec2}

Let $G$ be a non-universal profinite group with an infinite finitely generated simple quotient
$G/N$. Without loss of generality we may assume that $G$ is generated by
say $d$ elements. Let $G_1$ be the prosoluble
residual of $G$, i.e. the intersection of all open normal subgroups $N$
of $G$ such that
$G/N$ is soluble. Let $G_2$ be the the smallest closed normal subgroup
of $G_1$ such that $G_1/G_2$ is semisimple. (Alternatively $G_2$ is the intersection of all open
normal subgroups $M$ of $G_1$ such that $G_1/M$ is simple). Let 
$G_3=[G_2,_\omega G]$. This is a closed subgroup of $G$ by Theorem 1.4 of
\cite{NS} but strictly speaking we don't need this result here as we could
take $G_3$ as the closure of $[G_2,_\omega G]$ with the same effect. \medskip

The quotient $G/G_3$ is a pronilpotent-by-semisimple-by-prosoluble group
and therefore from Proposition \ref{main} and Theorem \ref{sol} we already know that Conjecture 1 holds for it. 

We have assumed that $N$ is a normal subgroup of $G$ such that
$G/N$ is infinite simple and finitely generated. Since $G/G_3$ doesn't have
infinite finitely generated images it is not possible that
$G_3 \leq N$ and therefore $NG_3=G$. Choose elements $n_1,n_2, \ldots , n_d$
from $N$ such that 

\[ \overline{\langle n_1, \ldots ,n_d \rangle} \ G_3 = G \]

Below we shall prove that there is a number $t$ depending on $d$ and $\alpha(G)$
such that

\begin{equation} \label{eq1}
G_3= \left([G,n_1] \cdot [G,n_2] \cdots [G,n_d] \right)^{*t}
\end{equation}

Since $[G,n_i] \subset N$ for each $i$ we see that $G_3 \leq N$ and hence
$G=NG_3=N$, contradiction.

It remains to prove equation (\ref{eq1}). Since both sides are closed subsets
of $G$ is is enough to show that (\ref{eq1}) is true in every finite image
of $G$. So let $\Gamma \geq \Gamma_2 \geq \Gamma_3$ be finite images of $G,G_2$ and $G_3$ respectively. It follows that $\Gamma_2$ is the semisimple
by soluble residual of $\Gamma$ and that $\Gamma_3= [G_2,_\omega G]$.

Recall the definition of acceptable subgroups in \cite{NS}, page 178:

A normal subgroup $N$ of $\Gamma$ is acceptable if
\medskip

1. $[N,G]=N$, and \medskip

2. If $A \leq B\leq N$ are two normal subgroups of $G$ inside $N$ then $B/A$
is not isomorphic to $S$ or $S \times S$ for a nonabelian finite simple group
$S$.

\bigskip

It follows that $\Gamma_3$ is an acceptable subgroup of $\Gamma$. The first
condition is clear. If the second condition fails then $\Gamma /C_\Gamma (A/B)$ is isomorphic to a subgroup of $\mathrm{Aut}(A/B)$ and for a simple group $S$ we have that $\mathrm{Aut}(S)$ and $\mathrm{Aut}(S\times S)= \mathrm{Aut}(S) \wr C_2$ are both semisimple-by-soluble. Therefore $C_\Gamma (A/B)$ contains $\Gamma_2
\geq \Gamma_3$ which implies that $A/B$ is abelian, contradiction.

The following result is proved as Key Theorem (C) in \cite{NS}. \bigskip

\textbf{KEY THEOREM (C)} \emph{Let $\Gamma$ be a $d$-generated finite group
and $H$ an acceptable normal subgroup of $G$. Suppose that $G=H \langle g_1,
\cdots ,g_r\rangle$. There is a number $t=t(d,\alpha(G))$ which depends only
on $d$ and $\alpha(G)$ such that}

\[ H= \left([H,g_1] [H,g_2] \cdots [H,g_r]\right)^{*t}. \]

\bigskip

Now Key Theorem (C) with $n_i$ in place of $g_i$ for $i=1,2,\ldots , d$ and $\Gamma_3$ in
place of $H$  gives 

\[ \Gamma_3=\left([\Gamma_3,n_1] [\Gamma_3,n_2] \cdots [\Gamma_3,n_d]\right)^{*t}. \]

Since this holds for any finite image $\Gamma$ of $G$ equation (\ref{eq1})
is proved and Theorem \ref{nou} follows. $\square$.

\section{Proofs of Theorems \ref{exdense} and \ref{g'}}\label{sec3}
We shall prove Theorem \ref{g'} first and then deduce Theorem \ref{exdense}.
We need the following result by Gasch\"{u}tz \cite{Ga}.

\begin{theorem} \label{g} Let $G$ be a profinite group generated topologically by $d$
elements. Suppose that $M$ is a closed normal subgroup of $G$ and $g_1,\ldots
,g_k$ are some $k \geq d$ elements of $G$ such that $\overline{\langle g_1,\ldots,g_k\rangle}M=G$.

Then there exist elements $g_i' \in g_iM$ such that $\overline{\langle g_1',\ldots,g_k'\rangle}=G$.
\end{theorem}

We shall use this result repeatedly in the following situation: suppose $N$ is a subgroup of $G$ containing elements
$g_1, \ldots, g_d$ whose images are dense in a quotient $G/M$ of $G$. Assume that for some normal closed subgroup $M_1$ of $G$ contained in $M$ we have $M \leq NM_1$. Then we can pick elements $g_{i}' \in N \cap g_iM$ such that 
$\overline{\langle g_1', \ldots, g_d' \rangle}M_1=G$.
\bigskip

\textbf{Proof of Theorem \ref{g'}.} Let $G_i$ for $i=1,2,3$ be the subgroups defined
in Section \ref{sec2}, namely $G_1$ is the intersection of all open
normal subgroups $V$ of $G$ such that $G/V$ is soluble, $G_2$ is the intersection
of the normal open subgroups $W$ of $G_1$ such that $G_1/W$ is simple, and
$G_3=[G_2,_\omega G]$. Let $K$ be the intersection of all open normal subgroups
of $G$ with $G/K$ simple and nonabelian. Thus $G_*=G' \cap K$. Put $G_s=G_1
\cap K= G_1 \cap G_*$ and let $G_v= [G_2,G]$. Then $G/G_s=G/G_1 \times G/K$
and the group $G_1/G_v$ is a quasisimple: a perfect central extension of
the semisimple group $G_1/G_2$. 
\medskip

\textbf{Step 1}: $G'G_1 \leq NG_1$.

By assumption we can find elements $g_1,\ldots, g_d$ of $N$ such that $\overline{\langle g_1,\ldots, g_d \rangle} G_*=G$. In addition
we assume that $d$ is more that the minimal number
of topological generators for $G$. By
Theorem \ref{sol} applied to $G/G_1$ and elements $g_iG_1$ it follows
that
\[ G'=  \left(\prod_{i=1}^d [G,g_{i}] \right)^{*(72d+46)}G_1 \subseteq NG_1,
\]
Theorem \ref{g} together with $G'/G_1 \leq NG_1/G_1$ implies that we can replace the elements
$g_i$ by some elements $g_{1,i} \in g_iG' \cap N$ such that $\overline{\langle g_{1,1}, \ldots, g_{1,d}\rangle} G_s=G$. 
\medskip

\textbf{Step 2:} $G_1 \leq NG_v$.

We can assume that $G_v=1$ in this step and will aim to show that $G_1 \leq N$. Now $G_1=G_1/G_v$ is a cover of the semisimple group $G_1/G_2$ and therefore a central product of finite quasismple groups. We can write $G_1 = P \circ Q$ as a central product of two qiasisimple groups $P$ and $Q$ where $Q= G_s'$ and $G_1 \cap K =G_s= \mathrm{Z}(G_1)Q$ and so $G/K \simeq G_1/G_s=P/ \mathrm{Z}(P)$.

First we shall show that $P \leq N$. The elements $\{g_{1,i}\}_{i=1}^d$ generate a dense subgroup of $G/G_s= G/G_1 \times G_1/G_s$. It follows that in every finite image of $\overline {G}$ of $G$ the images $\overline{g}_{1,i}$ generate a group which acts on the image $\overline{P}$ of $P$ as its full group of inner automorphisms.

We claim that in this case there is a number $c$ depending only on $\alpha(G)$ such that 
$\overline P = (\prod_{i=1}^d [\overline P, \overline g_{1,i}])^{*c}$ giving that
\begin{equation} P = (\prod_{i=1}^d [P,g_{1,i}])^{*c} \subseteq N.
\label{eq:1}
\end{equation}
\medskip

To see this first note that Theorem \ref{lie} of Liebeck and Shalev \cite{LS} has the following form in the special case when a quasisimple group has bounded Lie rank.

\begin{proposition}\label{basic} Let $S$ be a quasisimple group with $\alpha(S)=\alpha$. Then there is a number $c=c(\alpha) \in \mathbb N$ depending only on $\alpha$
such that $S=\mathcal C^{*c}$ for any noncentral conjugacy class $\mathcal C$ of $S$.
\end{proposition}

The claim (\ref{eq:1}) now follows from

\begin{proposition}\label{inner} Let $T$ be a finite central product of quasisimple groups and suppose $a_1,\ldots ,a_d \in \mathrm{Inn}(S)$ be $d$ generators of $\mathrm{Inn}(T)$. Then for $c=c(\alpha(T))$ as above we have
\[ T=\left(\prod_{i=1}^d[T,a_i]\right)^{*c}.\] 
\end{proposition}

Indeed for any quasisimple factor $S$ of $T$ there is at least one $a_j$ which acts as some nontrivial inner automorphism $s \mapsto s^g$ of $S$. Identifying $a_j$ with the noncentral element $g$ of $S$ we have $[S,a_j]=[S,g]=\mathcal C g$ where
$\mathcal C$ is the conjugacy class of $g^{-1}$ in $S$. Note that if $\mathcal
C_1$ and $\mathcal C_2$ are conjugacy classes of a group $S$ and $g \in S$
then 
\[g \mathcal C_1 =\mathcal C_1 g; \quad \mathcal C_1 \mathcal C_2= \mathcal
C_2 \mathcal C_1.\]

The result follows from Proposition \ref{basic}.

Therefore $P  \leq N$. 

The proof that $Q \leq N$ relies on the following result whose proof we defer
to the next section.

\begin{proposition}\label{qsimple} Let $l \in \mathbb{N}$ and let $T=S_1
S_2 \cdots S_m$ be a central product of quasisimple finite group with $\alpha(T)=\alpha$. Let $c=c(\alpha)$ be the number from Proposition
\ref{basic}. Then if $a_1,\ldots, a_d$ are automorphisms of $T$ such that $\cap_{i=1}^d C_T(a_i)$ does not contain any $S_j$ then 
\[ T= \left(\prod_{i=1}^d ([T,a_i] [T,a_i^{-1}])\right)^{*c}. \]
\end{proposition}

This is easily seen to imply that

\begin{equation} \label{eq2} Q \leq \left(\prod_{i=1}^d ([Q,g_{1,i}] [Q,g_{1,i}^{-1}])\right)^{*c} \subseteq N \end{equation}

Let $S$ be a quasisimple factor of $Q$. Suppose that all the elements $g_{1,i}$ centralized $S$. Since centralizers in profinite groups are closed subgroups
and $\overline{\langle g_{1,1}, \ldots g_{1,d} \rangle}G_s=G$ this gives that $C_G(S) S=G$. Therefore $S/\mathrm{Z}(S) \simeq G/C_G(S)$, hence $C_G(S) \geq G_* \geq S$, contradiction. Therefore at least one of the elements $g_{1,i}$ does not centralize $S$ and (\ref{eq2}) follows from Proposition \ref{qsimple} as before by passing to the finite images of $G$.

We showed $G_1/G_v=PQ \leq NG_v/G_v$ i.e. $G_1 \leq NG_v$.
Now find elements $g_{2,i} \in N \cap g_{1,i}G_1$ such that
$\overline{\langle \{g_{2,i}\}_{i=1}^d\rangle}G_v=G$.

\medskip

\textbf{Step 3:} $G_v \leq NG_3$.

By definition $[G_2,_\omega G]=G_3$ and by our choices of $g_{2,i}$ we have $\overline{\langle \{g_{2,i}\}_{i=1}^d
\rangle} G'=G$. Therefore by the profinite version of Lemma 2.4 of \cite{NS} with $G/G_3$ in place of
$G$ and $G_2/G_3$ in place of $H$ we have
\[ G_v=[G_2,G] \leq \left( \prod_{i=1}^d [G_2,g_{2,i}]\right) G_3  \subseteq NG_3.\]

Therefore we can replace $g_{2,i}$ by elements $g_{3,i} \in N \cap g_{2,i}G_v$ such that $\overline{\langle \{g_{3,i}\}_{i=1}^d \rangle} G_3=G$.
\medskip

\textbf{Step 4:} $G_3 \leq N$.

Now we are in position to apply Key Theorem (C): in every finite image $\overline{G}$ of $G$ the group $G_3$ maps onto an acceptable normal subgroup $\overline{G_3}$ of $\overline{G}$ and the images of the elements $g_{3,i}$ generate $\overline{G}/\overline{G_3}$. Conclusion:
\[ G_3= [G_3,G] = \left( \prod_{i=1}^d [G_3,g_{3,i}]\right)^{*t} \subseteq N,\]
where $t=t(d, \alpha(G))$ is the number from Key Theorem (C).

Altogether by combining the results of results of Steps 1 to 4 we deduce that $G' \leq N$. Theorem \ref{g'} follows.
$\square$

\bigskip

\textbf{Proof of Theorem \ref{exdense}}
The argument at the start of Section \ref{den} shows that if $G$ has infinite
abelianization or maps onto infinitely many simple groups then it has a proper
dense normal subgroup. Conversely suppose that $G/G'$ is finite and $G$ maps
onto finitely many simple groups. Then $G_*$ is an open subgroup of $G$.

Suppose that $N$ is a dense normal subgroup of $G$. We will show that $N=G$. First since $G/G_*$ is finite it follows that $N$ has a finitely generated
subgroup $D$ such that $DG_*=G$. Therefore we can apply Theorem \ref{exdense}
and so $N \geq G'$. But $G'$ is open hence so is $N$. Therefore $N=G$.
$\square$.

\subsection{Proof of Proposition \ref{qsimple}}
 \begin{lemma}\label{useful} Let $f \in \mathrm{Aut}(G)$. Then
 \[a^G =\{ g^{-1}ag \ | \ g \in G \}  \subseteq [G,f][G,f^{-1}]\]
 for all $a \in [G,f^{-1}]$.
 \end{lemma}
 \textbf{Proof:} For $x,y \in G$ put $z=xy^f$ and consider
 \[ [y,f] [z,f^{-1}]=y^{-1} x^{-1}x^{f^{-1}}y= [x,f^{-1}]^y.\]
 $\square$
 
 \medskip
 \textbf{Proof of Proposition \ref{qsimple}}
Every automorphism of $T$ permutes the subgroups $S_i$ among themselves. In
addition if $S$ is a quasisimple group
the action of $\mathrm{Aut}(S)$ on its simple factor $S/ \mathrm{Z}(S)$ is
faithful. Let $\pi_j$ be the projection $T \rightarrow T/C_T(S_j)= S_j/\mathrm{Z}(S_j)$.

From the condition on the automorphisms
$a_i$ it is easy to show that there exist elements $t_1,\ldots, t_d \in T$ with the following property: for any quasisimple subgroup $S_j$ of $T$, ($j\in \{1,2,\ldots ,m\}$) there exists
$i \in \{1,\ldots, d\}$ such that the projection $\pi_j(v_i)$ of the element $v_i:=[t_i,a_i^{-1}]$ on $S_j/\mathrm{Z}(S_j)$ is nontrivial. This is obvious if some $a_i$ leaves $S_j$ invariant
while acting nontrivially on it. Otherwise some $a_i$ will move $S_j$ to a different subgroup $S_l$ and then we can in fact achieve that the projections $\pi_l$ of $[t_i,a_i^{-1}]$ on
\textbf{each} of the groups $S_l/\mathrm{Z}(S_l)$ from the orbit of $S_j$ to be nontrivial. 

Now by Lemma \ref{useful} 
\[ \left( \prod_{i=1}^d v_i^T \right)^{*c} \subseteq \left(\prod_{i=1}^d ([T,a_i] [T,a_i^{-1}])\right)^{*c}.\] 
and on the other hand from Proposition \ref{lie} and the choice of $v_i=[t_i,a_i^{-1}]$
above we obtain 
\[ S_j \leq \left( \prod_{i=1}^d v_i^T \right)^{*c}\]
for all $j=1,\ldots, m$. Proposition \ref{qsimple} follows. $\square$

\end{document}